\documentclass[]{article}
\usepackage{amsfonts}
\usepackage{amsmath,amssymb}
\newtheorem{thm}{Theorem}[section]
\newtheorem{cor}[thm]{Corollary}
\newtheorem{prop}[thm]{Proposition}

\newtheorem{lem}[thm]{Lemma}
\newtheorem{Def}[thm]{Definition}

\newtheorem{ex}{Example}[section]

\newcommand{\be}{\begin{equation}}
\newcommand{\ee}{\end{equation}}
\newcommand{\ben}{\begin{enumerate}}
\newcommand{\een}{\end{enumerate}}

\title{Homogeneous geodesics of left invariant Finsler metrics}
\author{%
  Dariush Latifi  \\
\\
  Department of Mathematics, Mohaghegh Ardabili University,
  \\
  P.O. Box. 56199-11367, Ardabil, Iran\\
    dlatifi@gmail.com }

\date{}

\begin{document}
\maketitle

\begin{abstract}
In this paper, we study the set of homogeneous geodesics of a
left-invariant Finsler metric on Lie groups. We first give a simple
criterion that characterizes geodesic vectors. As an application, we
study some geometric properties of bi-invariant Finsler metrics on
Lie groups. In particular a necessary and sufficient condition that
left-invariant Randers metrics are of Berwald type is given. Finally
a correspondence of homogeneous geodesics to critical points of
restricted Finsler metrics is given. Then results concerning the
existence homogeneous geodesics are obtained.
\end{abstract}
\\
\section{Introduction}

\

A classical problem of differential geometry is to study geodesics
of Riemannian manifolds $(M,g)$. Of particular interest are
geodesics with some special properties, for example homogeneous
geodesics. A geodesic of a Riemannian manifold $(M,g)$ is called
homogeneous if
it is an orbit of a one-parameter group of isometries of $M$. For results on homogeneous geodesics in
homogeneous Riemannian manifolds we refer to [8], [14], [12], [11].\\
Homogeneous geodesics have important applications to mechanics. For
example, the equation of motion of many systems of classical
mechanics reduces to the geodesic equation in an appropriate
Riemannian manifold $M$.\\ Geodesics of left-invariant Riemannian
metrics on Lie groups were studied by V. I. Arnold extending Euler's
theory of rigid-body motion [1]. A major part of V. I. Arnold's
paper is devoted to the study of homogeneous geodesics. Homogeneous
geodesics are called by V. I. Arnold "relative equilibriums ". The
description of such relative equilibria is important for qualitative
description of the behavior of the corresponding mechanical system
with symmetries. There is a big literature in mechanics devoted to
the investigation of relative equilibria. Homogeneous geodesics are
interesting also in pseudo-Riemannian geometry and light-like
homogeneous geodesics are of particular interest. For results on
homogeneous geodesics in homogeneous pseudo-Riemannian manifolds we
refer for example to [18], [19], [21], [7], [3], [6]. In [18], [21]
and [7], the authors study plane-wave limits(Penrose limits )of
homogeneous spacetimes along light-like homogeneous geodesics.\\
About the existence of homogeneous geodesics in a general
homogeneous Riemannian manifold, we have, at first, a result due to
V. V. Kajzer who proved that a Lie group endowed with a
left-invariant Riemannian metric admits at least one homogeneous
geodesic [10]. More recently O. Kowalski and J. Szenthe extended
this result to all homogeneous Riemannian manifolds [13]. An
extension of result of [13] to reductive homogeneous
pseudo-Riemannian manifolds has been also obtained [21], [6].
Homogeneous geodesics of left-invariant Lagrangians on Lie groups
were studied by J. Szenthe [24]. In this paper, we study the set of
homogeneous geodesics of a left-invariant Finsler metric on Lie
groups.

\section{Preliminaries}
\subsection{Finsler spaces}

In this section, we recall briefly some known facts about Finsler
spaces. For details, see [2], [23], [4]. \

Let $M$ be a n-dimensional $C^{\infty}$ manifold and
$TM=\bigcup_{x\in M}
 T_{x}M$ the tangent bundle. If the continuous  function $F:TM\longrightarrow
 R_{+}$ satisfies the conditions that it is $C^{\infty}$ on $TM\setminus
 \{0\}$; $F(tu)=tF(u)$ for all $t\geq 0$ and $u\in TM$, i.e, $F$ is
 positively homogeneous of degree one; and for any tangent vector $y\in T_{x}M\setminus \{0\}$, the following bilinear symmetric form $g_{y}:T_{x}M \times T_{x}M\longrightarrow R$ is positive definite :
$$g_{y}(u,v)=\frac{1}{2}\frac{\partial^{2}}{\partial s\partial t}[F^{2}(x,y+su+tv)]|_{s=t=0},$$
  then we say that $(M,F)$ is a Finsler manifold.\\
Let $$g_{ij}(x,y)=(\frac{1}{2}F^{2})_{y^{i}y^{j}}(x,y).$$ By the
homogeneity of $F$, we have
$$g_{y}(u,v)=g_{ij}(x,y)u^{i}v^{j},\hskip.5cm F(x,y)=\sqrt{g_{ij}(x,y)y^{i}y^{j}}.$$

Let $\gamma:[0,r]\longrightarrow M$ be a piecewise  $C^{\infty}$
curve. Its integral length is defined as
$$L(\gamma)=\int_{0}^{r}F(\gamma(t),\dot{\gamma}(t))dt.$$ For $x_{0},x_{1}\in
M$ denote by  $\Gamma(x_{0},x_{1})$  the set of all piecewise
$C^{\infty}$ curve $\gamma:[0,r]\longrightarrow M$ such that
$\gamma(0)=x_{0}$ and $\gamma(r)=x_{1}$. Define a map $d_{F}:M\times
M\longrightarrow [0,\infty)$ by \[d_{F}(x_{0},x_{1})=\inf_{\gamma\in
\Gamma(x_{0},x_{1})} L(\gamma).\] Of course, we have
$d_{F}(x_{0},x_{1})\geq 0$, where the equality holds if and only if
$x_{0}=x_{1}$; $d_{F}(x_{0},x_{2})\leq
d_{F}(x_{0},x_{1})+d_{F}(x_{1},x_{2})$. In general, since $F$ is
only a positive homogeneous function, $d_{F}(x_{0},x_{1})\neq
d_{F}(x_{1},x_{0})$, therefore $(M,d_{F})$ is only a non-reversible
metric space.

Let $\pi^{\ast}TM$ be the pull-back of the tangent bundle $TM$ by
$\pi: TM\setminus\{0\}\longrightarrow M$. Unlike the Levi-Civita
connection in Riemannian geometry, there is no unique natural
connection in the Finsler case. Among these connections on
$\pi^{\ast}TM$, we choose the \emph{Chern connection} whose
coefficients are denoted by $\Gamma^{i}_{jk}$(see[2,p.38]). This
connection is almost $g-$compatible and has no torsion. Here
$g(x,y)=g_{ij}(x,y)dx^{i}\otimes
dx^{j}=(\frac{1}{2}F^{2})_{y^{i}y^{j}}dx^{i}\otimes dx^{j}$ is the
Riemannian metric on the pulled-back bundle $\pi^{\ast}TM$.

The Chern connection defines the covariant derivative $D_{V}U$ of a
vector field $U\in \chi(M)$ in the direction $V\in T_{p}M$. Since,
in general, the Chern connection coefficients $\Gamma^{i}_{jk}$ in
natural coordinates have a directional dependence, we must say
explicitly that $D_{V}U$ is defined with a fixed reference vector.
In particular, let $\sigma :[0,r]\longrightarrow M$ be a smooth
curve with velocity field $T=T(t)=\dot{\sigma}(t)$. Suppose that $U$
and $W$ are vector fields defined along $\sigma$. We define $D_{T}U$
with \emph{reference vector} $W$ as
$$D_{T}U=\left[\frac{dU^{i}}{dt}+U^{j}T^{k}(\Gamma^{i}_{jk})_{(\sigma ,W)}\right]\frac{\partial}{\partial x^{i}}\mid_{\sigma(t)}.$$
A curve $\sigma:[0,r]\longrightarrow M$, with velocity
$T=\dot{\sigma}$ is a Finslerian \emph{geodesic} if\\

$D_{T}\left[\frac{T}{F(T)}\right]=0$ ,\hskip.5cm with reference vector $T$.\\

We assume that all our geodesics $\sigma (t)$ have been
parameterized to have constant Finslerian speed. That is, the length
$F(T)$ is constant. These geodesics are characterized by the
equation \\

$D_{T}T=0$ , \hskip.5cm  with reference vector $T$.\\ \\Since
$T=\frac{d\sigma^{i}}{dt}\frac{\partial}{\partial x^{i}}$, this
equation says that
$$\frac{d^{2}\sigma^{i}}{dt^{2}}+\frac{d\sigma^{j}}{dt}\frac{d\sigma^{k}}{dt}(\Gamma^{i}_{jk})_{(\sigma ,T)}=0.$$

If $U,V$ and $W$ are vector fields along a curve $\sigma$, which has
velocity $T=\dot{\sigma}$, we have the derivative rule
$$\frac{d}{dt}g_{_W}(U,V)=g_{_W}(D_{T}U,V)+g_{_W}(U,D_{T}V)$$
whenever $D_{T}U$ and $D_{T}V$ are with reference vector $W$ and one
of the following conditions holds:
\begin{description}
    \item[i)] U or V is proportional to W, or
    \item[ii)] W=T and $\sigma$ is a geodesic.
\end{description}

\subsection{Left-invariant Finsler metrics on Lie groups}

Let $G$ be a connected Lie group with Lie algebra $\frak{g}=T_{e}G$.
We may identify the tangent bundle $TG$ with $G\times \frak{g}$ by
means of the diffeomorphism that sends $(g,X)$ to
$(L_{g})_{\ast}X\in T_{g}G$.

\begin{Def}
A Finsler function $F:TG\longrightarrow R_{+}$ will be called
$G$-invariant if $F$ is constant on all $G$-orbits in $TG=G\times
\frak{g}$; that is, $F(g,X)=F(e,X)$ for all $g\in G$ and $X\in
\frak{g}$.
\end{Def}
The G-invariant Finsler functions on $TG$ may be identified with the
Minkowski norms on $\frak{g}$. If $F:TG\longrightarrow R_{+}$ is an
G-invariant Finsler function, then we may define
$\widetilde{F}:\frak{g}\longrightarrow R_{+}$ by
$\widetilde{F}(X)=F(e,X)$, where $e$ denotes the identity in $G$.
Conversely, if we are given a Minkowski norm
$\widetilde{F}:\frak{g}\longrightarrow R_{+}$, then $\widetilde{F}$
arises from an G-invariant Finsler function $F:TG\longrightarrow
R_{+}$ given by $F(g,X)=\widetilde{F}(X)$ for all $(g,X)\in G\times
\frak{g}$.\\
Let $G$ be a connected Lie group, $L:G\times G\longrightarrow G$ the
action being defined by the left-translations
$L_{g}:G\longrightarrow G$, $g\in G$ and $TL:G\times
TG\longrightarrow TG$ the action given by the tangent linear maps
$TL_{g}:TG\longrightarrow TG$, $g\in G$ of the left-translations. \\
A smooth vector field $X:TG-\{0\}\longrightarrow TTG$ is said to be
left-invariant if
$$TTL_{g}\circ X\circ TL^{-1}_{g}=X \hskip1cm \forall g \in G.$$
By a classical argument of calculus of variation we have the
following proposition.
\begin{prop}
If $F:TG\longrightarrow R_{+}$ is a left-invariant Finsler metric
then its geodesic spray $X$ is left-invariant as well.
\end{prop}

\section{Homogeneous geodesics of left invariant Finsler metrics}
\begin{Def}
Let $G$ be a connected Lie group, $\frak{g}=T_{e}G$ its Lie algebra
identified with the tangent space at the identity element,
$\widetilde{F}:\frak{g} \longrightarrow R_{+}$ a Minkowski norm and
$F$ the left-invariant Finsler metric induced by $\widetilde{F}$ on
$G$. A geodesic $\gamma :R_{+}\longrightarrow G$ is said to be
\emph{homogeneous} if there is a $Z\in \frak{g}$ such that
$\gamma(t)=exp(tZ)\gamma(0)$, $t\in R_{+}$ holds. A tangent vector
$X\in T_{e}G-\{0\}$ is said to be a \emph{geodesic vector} if the
1-parameter subgroup $t\longrightarrow exp(tX)$, $t\in R_{+}$, is a
geodesic of $F$.
\end{Def}
The geodesic defined by a geodesic vector is obviously a homogeneous
one. Conversely, let $\gamma$ be a geodesic with $\gamma(0)=g$ which
is homogeneous with respect to a 1-parameter group of
left-translations, namely $$\gamma(t)=exp(tY)g,\hskip.5cm t\in
R_{+},$$ then a homogeneous geodesic $\widetilde{\gamma}$ is given
by
\begin{eqnarray*}
  \widetilde{\gamma}(t) &=& L^{-1}_{g}\circ \gamma(t)=L^{-1}_{g}\circ R_{g}\circ exp(tY)\\
   &=&exp(Ad(g^{-1})tY).e=exp(Ad(g^{-1})tY)\widetilde{\gamma}(0),
\end{eqnarray*}
which means that $X=Ad(g^{-1})Y$ is a geodesic vector.\\
For results on homogeneous geodesics in homogeneous Finsler
manifolds we refer to [16]. The basic formula characterizing
geodesic vector in the Finslerian case was derived in [16], Theorem
3.1. In the following theorem we present a new elementary proof of
this theorem for left invariant Finsler metrics on Lie groups.

\begin{thm}
Let $G$ be a connected Lie group with Lie algebra $\frak{g}$, and
let $F$ be a left-invariant Finsler metric on $G$. Then $X\in
\frak{g}-\{0\}$ is a geodesic vector if and only if
$$g_{X}(X,[X,Z])=0$$ holds for every $Z\in
\frak{g}$.
\end{thm}
Proof: Following the conventions of [9] a left-invariant vector
field associated to an element $X$ in $T_{e}G$ is denoted by
$\widetilde{X}:G\longrightarrow TG$; that is
$\widetilde{X}_{x}=L_{x\ast}X$. For any left invariant vector fields
$\widetilde{X}, \widetilde{Y}, \widetilde{Z}$ on $G$, we have
\begin{equation}
\widetilde{Y}g_{\widetilde{X}}(\widetilde{Z},\widetilde{X})=g_{\widetilde{X}}(D_{\widetilde{Y}}\widetilde{Z},\widetilde{X})+g_{\widetilde{X}}(\widetilde{Z},D_{\widetilde{Y}}\widetilde{X})
\qquad \mbox{with reference $\widetilde{X}$}
\end{equation}
Similarly,
\begin{equation}
\widetilde{Z}g_{\widetilde{X}}(\widetilde{Y},\widetilde{X})=g_{\widetilde{X}}(D_{\widetilde{Z}}\widetilde{Y},\widetilde{X})+
g_{\widetilde{X}}(\widetilde{Y},D_{\widetilde{Z}},\widetilde{X})
\end{equation}
\begin{equation}
\widetilde{X}g_{\widetilde{X}}(\widetilde{Z},\widetilde{X})=g_{\widetilde{X}}(D_{\widetilde{X}}\widetilde{Z},\widetilde{X})+
g_{\widetilde{X}}(\widetilde{Z},D_{\widetilde{X}},\widetilde{X})
\end{equation}
All covariant derivatives have $\widetilde{X}$ as reference
vector.\\
Subtracting (2) from the summation of (1) and (3) we get
\begin{eqnarray*}
  g_{\widetilde{X}}(\widetilde{Z},D_{\widetilde{X}+\widetilde{Y}}\widetilde{X})+ g_{\widetilde{X}}(\widetilde{X}-\widetilde{Y},D_{\widetilde{Z}}\widetilde{X})&=& \widetilde{Y}g_{\widetilde{X}}(\widetilde{Z},\widetilde{X})-\widetilde{Z}g_{\widetilde{X}}(\widetilde{Y},\widetilde{X})+\widetilde{X}g_{\widetilde{X}}(\widetilde{Z},\widetilde{X})\\
   &&-g_{\widetilde{X}}([\widetilde{Y},\widetilde{Z}],\widetilde{X})-g_{\widetilde{X}}([\widetilde{X},\widetilde{Z}],\widetilde{X}),
\end{eqnarray*}
where we have used the symmetry of the connection, i.e.,
$D_{\widetilde{Z}}\widetilde{X}-D_{\widetilde{X}}\widetilde{Z}=[\widetilde{Z},\widetilde{X}]$.
Set $\widetilde{Y}=\widetilde{X}-\widetilde{Z}$ in the above
equation, we obtain
\begin{equation}
2g_{\widetilde{X}}(\widetilde{Z},D_{\widetilde{X}}\widetilde{X})=2\widetilde{X}g_{\widetilde{X}}(\widetilde{Z},\widetilde{X})-\widetilde{Z}g_{\widetilde{X}}(\widetilde{X},\widetilde{X})-2g_{\widetilde{X}}([\widetilde{X},\widetilde{Z}],\widetilde{X}).
\end{equation}
Since $F$ is left-invariant, $dL_{x}$ is a linear isometry between
the spaces $T_{e}G=\frak{g}$ and $T_{x}G$, $\forall x\in G$.
Therefore for any left-invariant vector field $\widetilde{X},
\widetilde{Z}$ on $G$, we have
$$g_{\widetilde{X}}(\widetilde{Z},\widetilde{X})=g_{X}(Z,X)$$
i.e., the functions $g_{\widetilde{X}}(\widetilde{Z},\widetilde{X})$
, $g_{\widetilde{X}}(\widetilde{X},\widetilde{X})$ are constant.
Therefore from (4) the following is obtained
$$g_{\widetilde{X}}(\widetilde{Z},D_{\widetilde{X}}\widetilde{X})\mid_{e}=-g_{\widetilde{X}}([\widetilde{X},\widetilde{Z}],\widetilde{X})\mid_{e}=-g_{X}([X,Z],X).$$
Consequently the assertion of the theorem follows.$\Box$ \\\\
The following Proposition is well known for left-invariant
Riemannian metrics.
\begin{prop}
Let $G$ be a connected Lie group furnished with a left-invariant
Finsler metric $F$. Then the following are equivalent,
\begin{enumerate}
    \item $F$ is right-invariant, hence bi-invariant.
    \item $F$ is $Ad(G)-$invariant.
    \item$g_{Y}([X,U],V)+g_{Y}(U,[X,V])+2C_{Y}([X,Y],U,V)=0$,\hskip.4cm$\forall \hskip.2cm Y\in \frak{g}-\{0\}, X, U, V
    \in\frak{g}$, where $C_{y}$ is the Cartan
 tensor of $F$ at $Y$.\\

If the Finsler structure $F$ is absolutely homogeneous, then one
also has.

    \item The inversion map $g\longrightarrow g^{-1}$ is an
isometry
    of $G$.
\end{enumerate}
\end{prop}
Proof: The equivalence of the first two assertion is routine, and we
omit the details. The equivalence between (1) and (3) is a result of
S. Deng and Z. Hou [5]. If $F$ is absolutely homogeneous, one can
check quite easily that (4) is equivalent to (1).$\Box$
\begin{cor}
If $G$ is a Lie group endowed with a bi-invariant Finsler metric,
then the geodesics through the identity of $G$ are exactly
one-parameter subgroups.
\end{cor}
Proof: Since $F$ is bi-invariant, we have
$$g_{Y}([X,U],V)+g_{Y}(U,[X,V])+2C_{Y}([X,Y],U,V)=0$$$\forall \hskip.2cm Y\in \frak{g}-\{0\}, X, U, V
    \in\frak{g}$. It follows from the
    homogeneity of $F$ that $C_{Y}(Y,V,W)=0$. So we have
    $$g_{Y}([X,Y],Y)=0.$$ The result now follows from the
    Theorem 3.2.$\Box$\\\\
A connected Finsler space $(M,F)$ is said to be symmetric [15] if to
each $p\in M$ there is associated an isometry
$s_{p}:M\longrightarrow M$ which is \begin{description}
    \item[(i)] involutive ($s_{p}^{2}$ is the identity).
    \item[(ii)]has $p$ as an isolated fixed point, that is, there
    is a neighborhood $U$ of $p$ in which $p$ is the only fixed
    point of $s_{p}$.
\end{description}
$s_{p}$ is called the symmetry of the point $p$.
\begin{thm}
Suppose $G$ is a Lie group with a bi-invariant absolutely
homogeneous Finsler metric, then $G$ is a symmetric Finsler space.
\end{thm}
Proof: Consider the smooth mapping $f:G\longrightarrow G$,
$f:x\longrightarrow x^{-1}$. Then $f_{\ast}:T_{e}G\longrightarrow
T_{e}G$ maps a vector $\xi \in T_{e}G$ to $-\xi$; in particular,
$df$ is an isometry of $\frak{g}=T_{e}G$. Clearly,
$f=R_{g^{-1}}fL_{g^{-1}}$. Therefore, $df_{g}:T_{g}G\longrightarrow
T_{g^{-1}}G$ is an isometry for any $g\in G$.\\ Let
$s_{g}(x)=gx^{-1}g,\hskip.3cm g,x \in G$. The mapping $s_{g}$ is an
isometry, because $s_{g}=R_{g}fR_{g^{-1}}$. Thus, $s_{g}$ is an
isometry of $G$, obviously fixing the point $g$. Furthermore,
$s_{g}^{2}(x)=g(gx^{-1}g)^{-1}g=x.$ To show that $s_{g}$ is the
symmetry used in the definition of a symmetric Finsler space, it
suffices to show that $(s_{g})_{\ast}\xi=-\xi$ whenever $\xi \in
T_{g}G$. \\ Let us start with the case $g=e$. Let
$\xi=\frac{d}{dt}\gamma (t)|_{t=0}\in T_{e}G$ where $\gamma (t)$ is
a one-parameter subgroup of $G$. Then $\gamma (t)^{-1}=\gamma (-t)$,
and $(s_{e})_{\ast}\xi =\frac{d}{dt}\mid_{t=0}\gamma (-t)=-\xi$.
Now, if $\xi$ is in $T_{e}G$ for an arbitrary $g\in G$, then
$ds_{g}=dR_{g}dfdR_{g^{-1}}$, so
$ds_{g}(\xi)=dR_{g}(df(dR_{g^{-1}}(\xi)))=dR_{g}(-dR_{g^{-1}}(\xi))=-\xi$.\\$\Box$

\

Let $M$ be a smooth n-dimensional manifold, a Randers metric on $M$
consists of a Riemannian metric
$\widetilde{a}=\widetilde{a}_{ij}dx^{i}\otimes dx^{j}$ on $M$ and a
1-form $b=b_{i}dx^{i}$, [2], [22]. Here $\widetilde{a}$ and $b$
define a function $F$ on $TM$ by
$$F(x,y)=\alpha(x,y)+\beta(x,y) \hskip2cm x\in M ,y\in T_{x}M$$
where $\alpha(x,y)=\sqrt{\widetilde{a}_{ij}y^{i}y^{j}}$ ,
$\beta(x,y)=b_{i}(x)y^{i}$. $F$ is Finsler structure if
$\|b\|=\sqrt{b_{i}b^{i}}<1$ where $b^{i}=\widetilde{a}^{ij}b_{j}$,
and $(\widetilde{a}^{ij})$ is the inverse of $(\widetilde{a}_{ij})$.
The Riemannian metric $\widetilde{a}=\widetilde{a}_{ij}dx^{i}\otimes
dx^{j}$ induces the musical bijections between 1-forms and vector
fields on $M$, namely $\flat :T_{x}M\longrightarrow T_{x}^{\ast}M$
given by $y\longrightarrow \widetilde{a}_{x}(y,\circ)$ and its
inverse $\sharp :T_{x}^{\ast}M\longrightarrow T_{x}M$. In the local
coordinates we have $$(y^{b})_{i}=\widetilde{a}_{ij}y^{j}\hskip1cm
y\in
T_{x}M$$$$(\theta^{\sharp})^{i}=\widetilde{a}^{ij}\theta_{j}\hskip1cm
\theta\in T_{x}^{\ast}M$$ Now the corresponding vector field to the
1-form $b$ will be denoted by $b^{\sharp}$, obviously we have

$\|b\|=\|b^{\sharp}\|$ and
$\beta(x,y)=(b^{\sharp})^{\flat}(y)=\widetilde{a}_{x}(b^{\sharp},y)$.
Thus a Randers metric $F$ with Riemannian metric
$\widetilde{a}=\widetilde{a}_{ij}dx^{i}\otimes dx^{j}$ and 1-form
$b$ can be showed by
$$F(x,y)=\sqrt{\widetilde{a}_{x}(y,y)}+\widetilde{a}_{x}(b^{\sharp},y)\hskip1cm x\in M ,y\in T_{x}M$$
where $\widetilde{a}_{x}(b^{\sharp},b^{\sharp})<1 \hskip1cm \forall
x\in M$.
\\\\ Let
$F(x,y)=\sqrt{\widetilde{a}_{x}(y,y)}+\widetilde{a}_{x}(X,y)$ be a
left invariant Randers metric. It is easy to check that the
underlying Riemannian metric $\widetilde{a}$ and the vector field
$X$ are also left invariant.
\begin{thm}
Let $G$ be a Lie group with a left-invariant Randers metric $F$
defined by the Riemannian metric
$\widetilde{a}=\widetilde{a}_{ij}dx^{i}\otimes dx^{j}$ and the
vector field $X$. Then the Randers metric $F$ is of Berwald type if
and only if $ad_{X}$ is skew-adjoint with respect to $\widetilde{a}$
and $\widetilde{a}(X,[\frak{g},\frak{g}])=0$.
\end{thm}
Proof: For all $Y,Z\in \frak{g}$,
\begin{equation}
2\widetilde{a}(Y,\nabla_{Z}X)=\widetilde{a}(Z,[Y,X])+\widetilde{a}(X,[Y,Z])-\widetilde{a}(Y,[X,Z]).
\end{equation}
where $\nabla$ is the Levi-Civita connection of
$(M,\widetilde{a})$.\\
 If $ad_{X}$ is skew-adjoint, then first and last terms of
(5) sum to $0$. If additionally
$\widetilde{a}(X,[\frak{g},\frak{g}])=0$, then the middle term is
also $0$. So $\nabla_{Z}X=0$ for all $Z\in \frak{g}$, which means
that $X$ is parallel. By theorem 11.5.1. of [2] the Randers metric
is of Berwald type if
and only if $X$ is parallel with respect to $\widetilde{a}$.\\
Conversely, assume that the Randers metric is of Berwald type, so
the left side of (5) equals $0$ for all $Y,Z \in \frak{g}$. When
$Y=Z$, this yields $2\widetilde{a}(Y,[Y,X])=0$ for all $Y\in
\frak{g}$, which implies that $ad_{X}$ is skew-adjoint. This
property makes the first and third terms of (5) sum to zero, so
$\widetilde{a}(X,[Y,Z])=0$ for all $Y,Z \in \frak{g}$. In other
words
$\widetilde{a}(X,[\frak{g},\frak{g}])=0$.$\Box$\\\\
By a simple modification of the previous procedure, we can easily
obtain the following.
\begin{thm}
Let $(M=\frac{G}{H},F)$ be a homogeneous Randers space with $F$
defined by the Riemannian metric $\widetilde{a}$ and the vector
field $X$. Let $\frak{m}$ be the orthogonal complement of $\frak{h}$
in $\frak{g}$ with respect to the inner product induced on
$\frak{g}$ by $\widetilde{a}$. Then the Randers metric $F$ is of
Berwald type if and only if $\left(ad_{X}\right)_{\frak{m}}$ is
skew-adjoint  and
$\widetilde{a}\left(X,[\frak{m},\frak{m}]_{\frak{m}}\right)=0$,
where $\left(ad_{X}\right)_{\frak{m}}$ denotes
$\left(ad_{X}\right)_{\frak{m}}:\frak{m}\longrightarrow \frak{m}$,
$\left(ad_{X}\right)_{\frak{m}}(y)=[X,y]_{\frak{m}}$.

\end{thm}
\begin{thm}
Let $G$ be a Lie group with a bi-invariant Randers metric $F$
defined by the Riemannian metric
$\widetilde{a}=\widetilde{a}_{ij}dx^{i}\otimes dx^{j}$ and the
vector field $X$. Then the Randers metric $F$ is of Berwald type.
\end{thm}
Proof:  Let
$F(p,y)=\sqrt{\widetilde{a}_{p}(y,y)}+\widetilde{a}_{p}(X,y)$.\\ Now
for $s,t \in R$
\begin{eqnarray*}
  F^{2}(y+su+tv) &=& \widetilde{a}(y+su+tv,y+su+tv)+\widetilde{a}^{2}(X,y+su+tv) \\
   & & +2\sqrt{\widetilde{a}(y+su+tv,y+su+tv)} \widetilde{a}(X,y+su+tv)
\end{eqnarray*}
By definition
$$g_{y}(u,v)=\frac{1}{2}\frac{\partial^{2}}{\partial r \partial
s}F^{2}(y+ru+sv)\mid_{r=s=0}.$$ So by a direct computation we get

\begin{eqnarray*}
  g_{y}(u,v) &=& \widetilde{a}(u,v)+
  \widetilde{a}(X,u)\widetilde{a}(X,v)\\\\
  && +\frac{\widetilde{a}(u,v)\widetilde{a}(X,y)}{\sqrt{\widetilde{a}(y,y)}}-\frac{\widetilde{a}(v,y)\widetilde{a}(u,y)\widetilde{a}(X,y)}{\widetilde{a}(y,y)\sqrt{\widetilde{a}(y,y)}}
  \\\\
   &&+\frac{\widetilde{a}(X,v)\widetilde{a}(u,y)}{\sqrt{\widetilde{a}(y,y)}}+\frac{\widetilde{a}(X,u)\widetilde{a}(v,y)}{\sqrt{\widetilde{a}(y,y)}}.\\
\end{eqnarray*}
So for all $y,z \in \frak{g}$ we have
\begin{eqnarray*}
  g_{y}(y,[y,z]) &=& \widetilde{a}(y,[y,z])+\widetilde{a}(X,y)\widetilde{a}(X,[y,z])
  \\ \\
   &&
   +\frac{\widetilde{a}(y,[y,z])\widetilde{a}(X,y)}{\sqrt{\widetilde{a}(y,y)}}+\widetilde{a}(X,[y,z])\sqrt{\widetilde{a}(y,y)}\\\\
   &=&\widetilde{a}(y,[y,z])\left(1+\frac{\widetilde{a}(X,y)}{\sqrt{\widetilde{a}(y,y)}}\right)\\\\
   &&+\widetilde{a}(X,[y,z])\left(\widetilde{a}(X,y)+\sqrt{\widetilde{a}(y,y)}\right).
\end{eqnarray*}
So we have

\begin{equation}
  g_{y}(y,[y,z]) = \widetilde{a}(y,[y,z])\left(\frac{F(y)}{\sqrt{\widetilde{a}(y,y)}}\right)+\widetilde{a}(X,[y,z])F(y)
\end{equation}
Since $\widetilde{a}$ is bi-invariant, $\widetilde{a}(y,[y,z])=0$
and $ad(x)$ is skew-adjoint for every $x\in \frak{g}$. Since $F$ is
bi-invariant, $g_{y}(y,[y,z])=0$. So From (6) we get
$\widetilde{a}(X,[y,z])=0$ for all $y,z \in \frak{g}$. Therefore, by
Theorem 3.6, we see that $(G,F)$ is of Berwald type.$\Box$

\begin{cor}
Let $G$ be a Lie group with a left-invariant Randers metric $F$
defined by the Riemannian metric
$\widetilde{a}=\widetilde{a}_{ij}dx^{i}\otimes dx^{j}$ and the
vector field $X$. If the Randers metric $F$ is of Berwald type then
$X$ is a geodesic vector.
\end{cor}
The following lemma can be found in [20, p.301].
\begin{lem}\emph{(Milnor)}
Let $G$ be a Lie group endowed with a left-invariant Riemannian
metric $\widetilde{a}$. If $x \in \frak{g}$ is
$\widetilde{a}-$orthogonal to the commutator ideal $[\frak{g},
\frak{g}]$, then $Ricci(x)\leq 0$, with equality if and only if
$ad_{x}$ is skew-adjoint with respect to $\widetilde{a}$.
\end{lem}
\begin{cor}
Let $G$ be a Lie group with a left-invariant Randers metric $F$
defined by the Riemannian metric
$\widetilde{a}=\widetilde{a}_{ij}dx_{i}\otimes dx_{j}$ and the
vector field $X$. If the Randers metric $F$ is of Berwald type then
the Ricci curvature of $\widetilde{a}$ in the direction
$u=\frac{X}{\sqrt{\widetilde{a}(X,X)}}$ is zero.
\end{cor}
Proof: The corollary is a direct consequence of Theorem 3.6 and
Lemma 3.10 .\\$\Box$
\section{Homogeneous geodesics and the critical points of the restricted Finsler function}
Let $G$ be a connected Lie group, $\frak{g}=T_{e}G$ its Lie algebra,
$Ad:G\times \frak{g} \longrightarrow \frak{g}$ the adjoint action,
$G(X)=\{Ad(g)X\mid g\in G\}\subset \frak{g}$ the orbit of an element
$X\in \frak{g}$ and $G_{X}<G$ the isometry subgroup at $X$. The set
$\frac{G}{G_{X}}$ of left-cosets of $G_{X}$ endowed with its
canonical smooth manifold structure admits the canonical left-action
$$\Lambda:G\times \frac{G}{G_{X}}\longrightarrow \frac{G}{G_{X}}
\hskip1cm (g,aG_{X})\longrightarrow gaG_{X},$$ which is also smooth.
Moreover, a smooth bijection $\rho:\frac{G}{G_{X}}\longrightarrow
G(X)$ is defined by $\rho(aG_{X})=Ad(a)X$ which thus yields an
injective immersion into $\frak{g}$ which is equivariant
with respect to the actions $\Lambda$ and $Ad$.\\
Now consider a Minkowski norm $\widetilde{F}:\frak{g}\longrightarrow
R$, then $F$ defines a left-invariant Finsler metric on $G$ by
$$F(x,U)=\widetilde{F}(dL_{x^{-1}}U), \hskip.5cm U\in T_{x}G,$$ where $L_{x}:G\longrightarrow
G$ is the left translation by $x\in G$. Let
$Q(Z)=\widetilde{F}^{2}(Z)$, \hskip.2cm $Z\in \frak{g}$. Using the
formula $\widetilde{F}(Z)=\sqrt{g_{Z}(Z,Z)}$, we have $Q(Z)=g_{Z}(Z,Z)$.\\
The smooth function $q=Q\circ\rho:\frac{G}{G_{X}}\longrightarrow R$
will be called the restricted Minkowski norm on $\frac{G}{G_{X}}$.\\
In the following, we give an extension of results of [25] to
left-invariant Finsler metrics. We use some ideas from [25], [26] in
our proofs.
\begin{thm}
Let $G$ be a connected Lie group and $\widetilde{F}$ a Minkowski
norm on its Lie algebra $\frak{g}$. For $X\in \frak{g}-\{0\}$ let
$U\in \frak{g}$ be such that $X\in G(U)$ for the corresponding
adjoint orbit and let $gG_{U}\in \frac{G}{G_{U}}$ be the unique
coset with $\rho(gG_{U})=X$. Then $X$ is a geodesic vector if and
only if $gG_{U}$ is a critical point of $q=Q\circ \rho$ the
restricted Minkowski norm on $\frac{G}{G_{U}}$.
\end{thm}
Proof: The coset $gG_{U}$ is a critical point of $q$ if and only if
$vq=0$ for $v\in T_{gG_{U}}(\frac{G}{G_{U}})$. But as
$\frac{G}{G_{U}}$ is homogeneous, for each $v$ there is a $Z\in
\frak{g}$ such that $v=\widetilde{Z}(gG_{U})$ where
$\widetilde{Z}:\frac{G}{G_{U}}\longrightarrow T(\frac{G}{G_{U}})$ is
the infinitesimal generator of the action $\Lambda$ corresponding to
$Z$. Consider also the infinitesimal generator $\widehat{Z}:
\frak{g}\longrightarrow T\frak{g}$ of the adjoint action
corresponding to $Z$. Since the injective immersion $\rho$ is
equivariant with respect to the action $\Lambda$ and $Ad$ the
following holds: $\widehat{Z}\circ\rho=T\rho\circ\widetilde{Z}$. But
then the following is valid:
\begin{eqnarray*}
  v(q) &=& \widetilde{Z}(q)\mid_{gG_{U}}=\widetilde{Z}(Q\circ p)\mid_{gG_{U}}
  \\\\
   &=& \left(T\rho\widetilde{Z}\right)\mid_{gG_{U}}Q=(\widehat{Z}\circ\rho)\mid_{gG_{U}}Q \\
   &=& \left(\frac{d}{dt}\mid_{t=0}(Ad(exptZ)X)\right)Q \\
   &=& \frac{d}{dt}\mid_{t=0}Q(Ad(exptZ)X) \\
   &=&
   \frac{d}{dt}\mid_{t=0}g_{_{Ad(exptZ)X}}(Ad(exptZ)X,Ad(exptZ)X)\\\\
   &=& g_{X}([Z,X],X)+g_{X}(X,[Z,X])+2C_{X}([Z,X],X,X) \\\\
   &=& 2g_{X}([Z,X],X),
\end{eqnarray*}
where $C_{X}$ is the Cartan tensor of $F$ at $X$. It follows from
the homogeneity of $F$ that $C_{X}([Z,X],X,X)=0$. Since the map
$\alpha:\frak{g}\longrightarrow T_{gG_{U}}(\frac{G}{G_{U}})$,
$Z\longrightarrow \widetilde{Z}(gG_{U})$ is an epimorphism, the
assertion of the theorem follows.$\Box$
\begin{cor}
Let $G$ be a compact connected semi-simple Lie group and
$\widetilde{F}$ a Minkowski norm on its Lie algebra $\frak{g}$. Then
each orbit of the adjoint action $Ad:G\times \frak{g}\longrightarrow
\frak{g}$ contains at least two geodesic vectors.
\end{cor}
Proof: Consider an orbit $G(X)$ of the adjoint action, the
corresponding coset manifold $\frac{G}{G_{X}}$ and the injective
immersion $\rho:\frac{G}{G_{X}}\longrightarrow \frak{g}$. Since $G$
is compact and semi-simple then the manifold $\frac{G}{G_{X}}$
becomes compact, and the restricted Minkowski norm
$q=Q\circ\rho:\frac{G}{G_{X}}\longrightarrow R$ has at least two
critical points.$\Box$ \\\\
The following corollary is a consequence of the preceding corollary.
Two geodesics are considered different if their images are
different.
\begin{cor}
Let $G$ be compact connected semi-simple Lie group of
\hskip.1cm$rank\geq2$ and $\widetilde{F}$ a Minkowski norm on its
Lie algebra. Then the left-invariant Finsler metric $F$ induced by
$\widetilde{F}$ on $G$ has infinitely many homogeneous geodesic
issuing from the identity element.
\end{cor}
Proof: The proof is similar to the Riemannian case, so we omit it
[25].$\Box$
\\\\
\section{Some examples}

\begin{ex}
\end{ex}
Let $G$ be a three-dimensional connected Lie group endowed with a
left-invariant Riemannian metric $\widetilde{a}$.
\begin{enumerate}
    \item Let $G$ be an unimodular Lie group. According to a result
    due to J. Milnor (see[20,Theorem 4.3, p.305],[17]) there exist an
    orthonormal basis $\{e_{1}, e_{2}, e_{3}\}$ of the Lie algebra
    $\frak{g}$ such that $$[e_{1}, e_{2}]=\lambda_{3}e_{3}, \hskip.5cm [e_{2}, e_{3}]=\lambda_{1}e_{1},\hskip.5cm[e_{3}, e_{1}]=\lambda_{2}e_{2}.$$
    Let $F$ be a left invariant Randers metric on $G$ defined by the
    Riemannian metric $\widetilde{a}$ and the vector field $X=\epsilon
    e_{1}$, $0 <\epsilon<1$ i.e.
    $$F(p,y)=\sqrt{\widetilde{a}_{p}(y,y)}+\widetilde{a}_{p}(X,y).$$
    We note, by using Theorem 3.6, that $(G,F)$ is not of the
    Berwald type. We want to describe all geodesic vectors of
    $(G,F)$.\\ For $s,t \in R$ \begin{eqnarray*}
  F^{2}(y+su+tv) &=& \widetilde{a}(y+su+tv,y+su+tv)+\widetilde{a}^{2}(X,y+su+tv) \\
   & & +2\sqrt{\widetilde{a}(y+su+tv,y+su+tv)} \widetilde{a}(X,y+su+tv)
\end{eqnarray*}
By definition
$$g_{y}(u,v)=\frac{1}{2}\frac{\partial^{2}}{\partial r \partial
s}F^{2}(y+ru+sv)\mid_{r=s=0}.$$ So by a direct computation we get
\begin{eqnarray*}
  g_{y}(u,v) &=& \widetilde{a}(u,v)+
  \widetilde{a}(X,u)\widetilde{a}(X,v)\\\\
  && +\frac{\widetilde{a}(u,v)\widetilde{a}(X,y)}{\sqrt{\widetilde{a}(y,y)}}-\frac{\widetilde{a}(v,y)\widetilde{a}(u,y)\widetilde{a}(X,y)}{\widetilde{a}(y,y)\sqrt{\widetilde{a}(y,y)}} \mbox{\qquad
  }\\\\
   &&+\frac{\widetilde{a}(X,v)\widetilde{a}(u,y)}{\sqrt{\widetilde{a}(y,y)}}+\frac{\widetilde{a}(X,u)\widetilde{a}(v,y)}{\sqrt{\widetilde{a}(y,y)}}.\\
\end{eqnarray*}
So for all $z\in \frak{g}$ we have
\begin{equation}
g_{y}(y,[y,z])=
\widetilde{a}\left(X+\frac{y}{\sqrt{\widetilde{a}(y,y)}} \hskip.1cm,
[y,z]\right)F(y)
\end{equation}
Using Theorem 3.2 and (7) we can check easily that $e_{1}$ is a
geodesic vector.\\ By using Theorem 3.2 and (7) a vector
$y=y_{1}e_{1}+y_{2}e_{2}+y_{3}e_{3}$ of $\frak{g}$ is a geodesic
vector if and only if $$\widetilde{a}\left(\epsilon
e_{1}+\frac{y_{1}e_{1}+y_{2}e_{2}+y_{3}e_{3}}{\sqrt{y_{1}^{2}+y_{2}^{2}+y_{3}^{2}}}
\hskip.1cm, [y_{1}e_{1}+y_{2}e_{2}+y_{3}e_{3} , e_{j}]\right)=0$$
for each $j=1, 2, 3.$\\ So we get :
$$(\lambda_{2}-\lambda_{3})y_{2}y_{3}=0,$$
$$-\epsilon y_{3}\lambda_{1}-\frac{1}{\sqrt{y_{1}^{2}+y_{2}^{2}+y_{3}^{2}}}y_{1}y_{3}\lambda_{1}+\frac{1}{\sqrt{y_{1}^{2}+y_{2}^{2}+y_{3}^{2}}}y_{1}y_{3}\lambda_{3}=0,$$
$$-\epsilon y_{2}\lambda_{1}+\frac{1}{\sqrt{y_{1}^{2}+y_{2}^{2}+y_{3}^{2}}}y_{1}y_{2}\lambda_{1}-\frac{1}{\sqrt{y_{1}^{2}+y_{2}^{2}+y_{3}^{2}}}y_{1}y_{2}\lambda_{2}=0.$$
As a special case, if $\lambda_{1}=\lambda_{2}=\lambda_{3}\neq 0$ we
conclude that all geodesic vectors $y$ are those from the set
$Span\{e_{1}\}$. Consequently, there is only one homogeneous
geodesic.
    \item Let $G$ be a non-unimodular Lie group. According to a
    result due to J. Milnor (see[20, Lemma 4.10, p.309],[17]) there
    exists an orthogonal basis $\{e_{1}, e_{2}, e_{3}\}$ of the Lie
    algebra $\frak{g}$ such that
    $$[e_{1},e_{2}]=\alpha e_{2}+\beta e_{3},\hskip.5cm[e_{2},e_{3}]=0,\hskip.5cm [e_{1},e_{3}]=\gamma e_{2}+\delta e_{3}, $$
    where $\alpha , \beta, \gamma, \delta $ are real numbers such
    that the matrix
$$
\left(%
\begin{array}{cc}
  \alpha & \beta \\
  \gamma & \delta \\
\end{array}%
\right)
$$
has trace $\alpha+\delta =2$ and $\alpha \gamma + \beta\delta =0.$
Let $F$ be a left invariant Randers metric on $G$ defined by the
Riemannian metric $\widetilde{a}$ and the vector field $X=\epsilon
e_{1}, \hskip.2cm 0<\epsilon <1.$\\
By using Theorem 3.2 and (7), a vector
$y=y_{1}e_{1}+y_{2}e_{2}+y_{3}e_{3}$ of $\frak{g}$ is a geodesic
vector if and only if $$\widetilde{a}\left(\epsilon
e_{1}+\frac{y_{1}e_{1}+y_{2}e_{2}+y_{3}e_{3}}{\sqrt{y_{1}^{2}+y_{2}^{2}+y_{3}^{2}}}
\hskip.1cm, [y_{1}e_{1}+y_{2}e_{2}+y_{3}e_{3} , e_{j}]\right)=0$$
for each $j=1, 2, 3.$\\
This condition leads to the system of equations\\\\
$$y_{2}(-\alpha y_{2}-\gamma y_{3})+y_{3}(-y_{2}\beta -\delta y_{3})=0,$$
$$y_{1}y_{2}\alpha+y_{1}y_{3}\beta=0,$$
$$y_{1}y_{2}\gamma+y_{1}y_{3}\delta=0.$$
Putting $\alpha=2, \delta=0, \gamma=0$ the above equations take the
form $$2y_{2}\left(y_{2}+\frac{\beta}{2}y_{3}\right)=0,$$
$$2y_{1}\left(y_{2}+\frac{\beta}{2}y_{3}\right)=0.$$
So a vector $y$ of $\frak{g}$ is a geodesic vector if and only if :
\begin{description}
    \item[-] $y\in Span(e_{1}, e_{3})$ for $\beta=0.$
    \item[-] $y\in Span(e_{1})\bigcup Span(e_{3})\bigcup
    Span(\frac{\beta}{2}e_{2}-e_{3})$ for $\beta \neq 0$
\end{description}
\end{enumerate}

\noindent

\end{document}